\documentclass[graybox]{svmult}

\usepackage{type1cm}        
\usepackage{makeidx}         
\usepackage{graphicx}        
\usepackage{multicol}        
\usepackage[bottom]{footmisc}

\usepackage{newtxtext}       %
\usepackage[varvw]{newtxmath}       

\usepackage[T1]{fontenc}
\usepackage[utf8]{inputenc}
\usepackage{amsmath}
\usepackage{pgfplots}
\usepackage{pgfplotstable}
\usepackage{cite}
\usepackage[hidelinks]{hyperref}

\usepackage{cleveref}
\usepackage{float}
\usepackage{lipsum}
\usepackage{algorithm}
\usepackage{algpseudocode}
\usepackage{tikz}
\usepackage{booktabs}
\usepackage{subcaption}
\usetikzlibrary{shapes.geometric, arrows.meta, positioning, fit, backgrounds, calc}

\pgfplotsset{compat=1.18}

\definecolor{color1}{RGB}{228,26,28}
\definecolor{color2}{RGB}{55,126,184}
\definecolor{color3}{RGB}{77,175,74}
\definecolor{color4}{RGB}{152,78,163}
\definecolor{color5}{RGB}{255,127,0}

\newcommand{\lo}{\text{lo}}
\newcommand{\up}{\text{up}}
\newcommand\figref[1]{\text{Figure \ref{#1}}}
\newcommand\tabref[1]{\text{Table \ref{#1}}}

\pgfplotsset{
	table/columns/Graph/.style={string type},
	filter graph/.style={
		x filter/.code={
			\pgfplotstablegetelem{\coordindex}{Graph}\of{\datatable}
			\edef\tempa{\pgfplotsretval}
			\edef\tempb{#1}
			\ifx\tempa\tempb
			\else
			
			\fi
		}
	}
}

\pgfplotsset{
	table/columns/h/.style={string type}, 
	filter h/.style={
		x filter/.code={
			\pgfplotstablegetelem{\coordindex}{h}\of{\datatable}
			\edef\tempa{\pgfplotsretval}
			\edef\tempb{#1}
			\ifx\tempa\tempb
			\else
			
			\fi
		}
	}
}

\makeindex             
\begin{document}
	\title*{Efficient generation of Gaussian random fields on metric graphs via domain decomposition and mass matrix lumping}
	\titlerunning{Efficient Simulation of Gaussian Random Fields on Metric Graphs}
	%
	\author{
		Mihály Kovács\inst{1,2,3} \and
		Gyula Molnár\inst{1} \and
		Máté András Száraz\inst{1}
	}
	\authorrunning{Kovács, Molnár, Száraz}
	\institute{
		\textsuperscript{1} Faculty of Information Technology and Bionics, P\'azm\'any P\'eter Catholic University, Pr\'ater utca 50/a., H-1083 Budapest, Hungary \\
		\textsuperscript{2} Department of Analysis and Operations Research, Budapest University of Technology and Economics, M\H{u}egyetem rkp. 3., H-1111 Budapest, Hungary\\
		\textsuperscript{3} Department of Mathematical Sciences,Chalmers University of Technology and University of Gothenburg, SE-41296 Gothenburg, Sweden
	}
	\maketitle
	\abstract{We consider Gaussian Random Fields (GRFs) on metric graphs defined implicitly as the stationary solution to a fractional SPDE driven by Gaussian white noise. Sampling from the finite element approximation requires the Cholesky factorization of the mass matrix, causing non-linear execution time explosions and massive memory fill-in on large graphs. Hence, we combine Neumann-Neumann graph decomposition with mass matrix lumping and demonstrate empirically, that our approach preserves exact theoretical convergence rates established in \cite{Bolin_2023} while achieving multi-order speedups and massive memory reductions.}
	\section{Introduction}
	Over the past several decades, differential operators on metric graphs have been increasingly utilized to model phenomena across a diverse array of scientific disciplines \cite{Alexander_1983,Bouchaud_1987,Flesia_1989,Kallianpur_1984,Cho_2018,Avdonin_2023,Avdonin_2019,Mehandiratta_2021,Stoll_2021}.
	While the SPDE approach rigorously defines Matérn fields on graphs, standard Finite Element approximations via the Balakrishnan integral \cite{Bolin_2023} face severe computational bottlenecks. Generating spatial white noise requires a Cholesky factorization of the consistent mass matrix, causing non-linear execution time explosions and massive memory fill-in on large graphs. We counteract these bottlenecks by combining Neumann-Neumann domain decomposition with mass matrix lumping, reducing noise generation to an $\mathcal{O}(N)$ operation with reduced memory bloat. Concurrently, domain decomposition condenses the global PDE to a Schur complement system on the topological vertices. This allows matrix-free resolution via a Preconditioned Conjugate Gradient (PCG) solver, while internal edge dynamics are resolved locally using the tridiagonal Thomas algorithm.	
	The paper is organized as follows: Section \ref{sec:preliminaries} overviews mathematical preliminaries with regards to Metric Graphs, the SPDE approach to generating Gaussian Random Fields, the Balakrishnan integral representation of a fractional inverse operators, finite element discretization, white noise generation in a FEM-context, mass matrix lumping, and finally, domain decomposition through taking the Schur complement system. Section \ref{sec:algorithmic} then details the algorithmic implementation and Section \ref{sec:numerical} demonstrates empirically that our approach preserves exact theoretical convergence rates while achieving multi-order speedups and massive memory reductions.
	\section{Preliminaries and Mathematical Setup} \label{sec:preliminaries}
	\subsection{Metric Graphs and Function Spaces}
	A compact metric graph $\Gamma = (V, E)$ consists of a finite set of vertices ${V}$ and edges $E$ \cite{Berkolaiko_2012, Arioli_2017}. Each edge $e \in E$ has length $\ell_e \in (0, \infty)$ and a local coordinate $x \in [0, \ell_e]$. A function $u$ on $\Gamma$ is represented as a vector of edge functions $u = [u_e]_{e\in E}$. We consider functions that are continuous on $\Gamma$ and satisfy the Neumann-Kirchhoff condition $\sum_{e\in E_v} u_e^{\prime}(v) = 0, \; v \in V$ \cite{Arioli_2017}. The relevant Hilbert space is defined as the direct sum over the edges $L^2(\Gamma) = \bigoplus_{e \in E} L^2(0,\ell_{e})$
	with the norm $\|u\|_{L^2(\Gamma)}^2 = \sum_{e\in E} \|u_e\|_{L^2(0,\ell_e)}^2$.
	\subsection{Gaussian Random Fields and the SPDE approach}
	A zero-mean Gaussian Random Field (GRF) $u$ on $\Gamma$ is a random field such that for $\{x_1, \dots, x_{k \in \mathbb{N}}\} \subset \Gamma$ we have that $[u(x_1), \dots, u(x_k)]^T$ is a zero mean Gaussian vector, the latter denoted by $ \mathcal{N}({0}, \boldsymbol{\Sigma})$, where $\boldsymbol{\Sigma}$ is the covariance matrix of the vector. One of the most popular family of GRFs in statistical modeling are Gaussian Mat\'ern fields. However, it is difficult to define valid isotropic covariance function in case of a metric the graph \cite{AMR}. In analogy to the case of $\mathbb{R}^n$ one may define a familiy of Gaussian random fields as the solution to the stochastic partial differential equation \cite{Lindgren_2011, Bolin_2023}:
	\begin{equation} \label{eq:fractional_spde}
		L^\beta u = (\kappa^2 - \Delta)^\beta u = \mathcal{W} \quad \text{on } \Gamma,
	\end{equation}
	where $\Delta$ is the Laplace operator acting on each edge, $\kappa > 0$ and $\beta \in (1/4,1)$ control main charactersitics of the field: the smoothness, variance and correlation length, and $\mathcal{W}$ is Gaussian white noise on $L^2(\Gamma)$. The solution to \eqref{eq:fractional_spde} can be the viewed as the generalization of Gaussian Mat\'ern fields to metric graphs \cite{Bolin_2023,Bolin_2024}.
	\subsection{The Balakrishnan Integral Representation}
	To circumvent computing $u = L^{-\beta} \mathcal{W}$ directly, we utilize the Balakrishnan integral representation of \eqref{eq:fractional_spde} \cite{Balakrishnan_1960}:
	\begin{equation} \label{eq:balakrishnan_integral}
		L^{-\beta} = \frac{\sin(\pi \beta)}{\pi} \int_{0}^{\infty} s^{-\beta} (L + sI)^{-1} \, ds.
	\end{equation}
	By applying the variable substitution $s = e^y$, and discretizing the integral with a uniform step size $k > 0$, the fractional solution is approximated by \cite{Bonito_2015}
	\begin{equation} \label{eq:balakrishnan_quadrature}
		u \approx \frac{k \sin(\pi \beta)}{\pi} \sum_{l=K^-}^{K^+} e^{(1-\beta) y_l} (L + e^{y_l} I)^{-1} \mathcal{W},
	\end{equation}
	where $y_l = l k$. 
	By utilizing the identity $w A^{-1} = (w^{-1} A)^{-1}$, we define the step-dependent scalar coefficients for the identity ($c_{I}^{(l)}$) and Laplacian ($c_{L}^{(l)}$) operators as:
	\begin{equation} \label{eq:coefficients}
		c_\text{const} = \frac{\pi}{2k\sin(\pi\beta)},
		\quad \quad c_{I}^{(l)} = c_\text{const} \cdot e^{-2\beta y_l},
		\quad \quad c_{L}^{(l)} = c_{I}^{(l)} e^{2y_l}.
	\end{equation}
	This allows us to rewrite the approximation as:
	$$ L^{-\beta} \approx \sum_{l=-K^-}^{K^+} \left( c_{I}^{(l)} I + c_{L}^{(l)} L \right)^{-1}. $$
	\subsection{Finite Element Discretization}
	For each quadrature step $l \in [-K^-, K^+]$, we must solve the continuous operator equation:
	$$ \left( c_{I}^{(l)} I + c_{L}^{(l)} L \right) u_l = \mathcal{W} $$
	which gets recast as a system of discrete linear algebraic equations
	\begin{equation} \label{eq:spde_approx_fem}
		{u} \approx \sum_{l=-K^-}^{K^+} {u}^{(l)} = \sum_{l=-K^-}^{K^+} {\underbrace{\left( c_{I}^{(l)} M + c_{L}^{(l)} G \right)}_{A^{(l)}}}^{-1} {W}
	\end{equation}
	with mass matrix $M$, stiffness matrix $G$, global system matrix $A^{(l)} \in \mathbb{R}^{N \times N}$, the statically generated discrete spatial white noise vector ${W}$. The numerical approximation of the fractional SPDE is then obtained by solving $A^{(l)} {u}^{(l)} = {W}$ for all $l$, and aggregating the resulting discrete vectors. Or, in block form:
	\begin{equation} \label{eq:spde_quadrature_block}
		u=
		\begin{bmatrix}
			u_I      \\
			u_\Gamma 
		\end{bmatrix} =
		\sum_{l} \left(
		\begin{bmatrix}
			A^{(l)}_{II}       & A^{(l)}_{I\Gamma}      \\
			A^{(l)}_{\Gamma I} & A^{(l)}_{\Gamma\Gamma} 
		\end{bmatrix}^{-1}
		\begin{bmatrix}
			W^{(l)}_I      \\
			W^{(l)}_\Gamma 
		\end{bmatrix}
		\right).
	\end{equation}
	\subsubsection*{Generation of the Discrete White Noise Vector}
	To simulate the fractional SPDE using the Finite Element Method, we project the continuous spatial white noise $\mathcal{W}(x)$ onto the basis functions $\{\phi_i\}_{i=1}^N$. The $i$-th component of the discrete noise vector is given by the $L^2$-projection $W_i = \int_\Gamma \mathcal{W}(x) \phi_i(x) dx$. Because $\mathcal{W}(x)$ is uncorrelated Gaussian white noise with spatial covariance $\delta(x-y)$, the covariance matrix of the projected vector $W$ corresponds exactly to the FEM consistent mass matrix, $\text{Cov}(W) = M$ \cite{Lindgren_2011, Bolin_2023}.	To generate a random vector with this exact geometric covariance structure computationally, we sample standard independent Gaussian noise $\boldsymbol{\xi} \sim \mathcal{N}({0}, I_N)$ and utilize the Cholesky decomposition of the symmetric positive-definite mass matrix, $M = R R^T$. The appropriately correlated discrete white noise vector is then obtained as:
	\begin{equation} \label{eq:wh_def}
		W = R \boldsymbol{\xi}.
	\end{equation}
	\subsection{Mass Lumping}
	The consistent mass matrix $M$ has entries $M_{ij} = \int_{\Gamma} \phi_i \phi_j \, dx$. Computing the Cholesky factor $R$ such that $M = RR^T$ for noise generation requires a global Cholesky factorization, which is computationally expensive \cite{Davis_2006}. We instead 
	use the diagonal lumped mass matrix: $M^\text{lumped}_{ii} = \sum_{j} M_{ij} = \int_{\Gamma} \phi_i \, dx$, and $M^\text{lumped}_{ij} = 0$ for $i \neq j$ \cite{Quarteroni_1994}. This allows the matrix square root $(M^\text{lumped})^{1/2}$ (the same as the Cholesky factor of $M^\text{lumped}$)  to be computed in $\mathcal{O}(N)$ operations via taking element-wise scalar square roots \cite{Golub_2013}.
	\subsection{Domain Decomposition}
	Solving the global block system \eqref{eq:spde_quadrature_block} directly is inefficient. Instead, non-overlapping domain decomposition methods reduce the global problem to an equation strictly on the interface (the topological vertices $\Gamma$) \cite{Toselli_2005}. By expressing the internal unknowns $u_I^{(l)}$ in terms of the boundary unknowns $u_\Gamma^{(l)}$ and substituting them into the second block row, we obtain the global Schur complement system:
	\begin{equation} \label{eq:global_schur}
		S^{(l)} u_\Gamma^{(l)} = g_\Gamma^{(l)},
	\end{equation}
	where $S^{(l)} = A_{\Gamma\Gamma}^{(l)} - A_{\Gamma I}^{(l)} (A_{II}^{(l)})^{-1} A_{I\Gamma}^{(l)}$ is the Schur complement matrix, and $g_\Gamma^{(l)} = W^{(l)}_\Gamma - A_{\Gamma I}^{(l)} (A_{II}^{(l)})^{-1} W^{(l)}_I$ is the condensed right-hand side.
	
	Because $S^{(l)}$ is symmetric positive-definite, \eqref{eq:global_schur} can be solved iteratively using the Preconditioned Conjugate Gradient (PCG) method \cite{Saad_2003, Toselli_2005}.
	\section{Algorithmic Implementation} \label{sec:algorithmic}
	\subsection{Phase 1: Initialization}
	After the arguments get parsed, and the graph gets loaded, the working mesh gets generated. Since we are discretizing a metric graph , every edge $e \in E$ has a fixed physical length $l_e$. For a target global minimum mesh resolution $h$ the number of uniform subintervals is calculated as $n_e = \left\lceil \frac{l_e}{h} \right\rceil$. Then, the local mesh size is $h_e = \frac{l_e}{n_e}$, and the total number of nodes (including boundary vertices) on edge $e$ is $N_e = n_e + 1$. With that, the total number of (topological-, and internal) vertices in the discretized graph is $N = |V| + \sum_{e \in E} (N_e - 2)$, where $V$ is the set of all topological graph vertices, and $E$ the set of all edges.
	\subsection{Phase 2: Global Setup and Noise Generation}
	Next, we generate the right hand side of \eqref{eq:spde_approx_fem}, that is, \eqref{eq:wh_def}. Herefore, we draw a sample from the standard Gaussian white noise $\boldsymbol{\xi} \sim \mathcal{N}({0}, I_N)$. Furthermore, we build the mass matrix $M = \left[M_{i,j}\right]$ as follows:
	\begin{enumerate}
		\item Consistent mass-matrix case:
		\begin{align}
			&M_{i,i} = \frac{2h_e}{3}, \quad e \in E, \;\; i \in I_e ,                                         \\
			&M_{i,j} = \frac{h_e}{6}, \quad e \in E, \;\; i,j \in I_e \text{ where } i,j \text{ are adjacent,} \\
			&M_{v,v} = \sum_{e \in E_v} \frac{h_e}{3}, \quad v \in V , \\
			&M_{v, i} = M_{i, v} = \frac{h_e}{6}, \quad v \in V, \;\; e \in E_v, \ \text{where } i = \text{adj}(v, e) , \\
			&M_{a,b} = 0 \quad \text{otherwise}.                                                                     
		\end{align}
		\item Lumped mass-matrix case:
		\begin{align}
			& M^{\text{lumped}}_{a,b} = \begin{cases}                                                                                                                                                                 
				\sum_{k} M_{a,k} & \text{if } a = b                                                                                                                                                                                        \\
				0                & \text{if } a \neq b                                                                                                                                                                                     
			\end{cases} \quad \text{, that is:}\\
			& M^{\text{lumped}}_{i,i} = M_{i,i} + \sum_j M_{i,j} = \frac{2h_e}{3} + \frac{h_e}{6} + \frac{h_e}{6} = h_e , \quad e \in E, \;\; i \in I_e , \\
			& M^{\text{lumped}}_{v,v} = M_{v,v} + \sum_j M_{v,j} = \left( \sum_{e \in E_v} \frac{h_e}{3} \right) + \left( \sum_{e \in E_v} \frac{h_e}{6} \right)=\sum_{e \in E_v} \frac{h_e}{2} , \quad v \in V \text{.}
		\end{align}
	\end{enumerate}
	Where $I_e$ the set of strictly internal nodes belonging to edge $e$, $E_v$ the set of edges incident to vertex $v \in V$. Furthermore, $\text{adj}(v, e)$ denotes the single internal node on edge $e$ that is immediately adjacent to vertex $v$.
	Then, we compute $R$ of \eqref{eq:wh_def} as the Cholesky-decomposition of $M$. In the lumped case, since $M^{\text{lumped}}$ is diagonal, its Cholesky-decomposition is obtained by taking element-wise square roots of the nonzero elements $\left[ \sqrt{M_{i,j}} \right]$ i.e. along the main diagonal. Finally, evaluating \eqref{eq:wh_def} is a matrix-vector multiplication, resulting in $W$, the right-hand side of \eqref{eq:spde_quadrature_block}. That is, $W = [{W}_I ~ {W}_\Gamma]^T = [{W}_{I,1} \dots {W}_{I,|E|} ~ {W}_\Gamma]^T$. As the very last step in this phase, we indeed partition and distribute $W$ into internal edge contributions ${W}_{I,e \in E}$, and the boundary node contribution ${W}_\Gamma$.
	Furthermore, the global stiffness matrix $G = \left[G_{i,j}\right]$ is assembled analogously from the standard 1D linear finite element gradients. Note that unlike the mass matrix, the stiffness matrix is evaluated exactly and is not subject to lumping:
	
	\begin{align}
		&G_{i,i} = \frac{2}{h_e} , \quad e \in E, \;\; i \in I_e,                                            \\
		&G_{i,j} = -\frac{1}{h_e}, \quad e \in E, \;\; i,j \in I_e, \text{ where } i,j \text{ are adjacent}, \\
		&G_{v,v} = \sum_{e \in E_v} \frac{1}{h_e}, \quad v \in V,                                                \\
		&G_{v, i} = G_{i, v} = -\frac{1}{h_e}, \quad v \in V, \;\; e \in E_v, \ \text{where } i = \text{adj}(v, e), \\
		&G_{a,b} = 0 \quad \text{otherwise.}                                                                            
	\end{align}
	\subsection{Phase 3: The Outside Loop: Balakrishnan Quadrature}
	\subsubsection*{Loop setup}
	Next, we prepare for evaluating \eqref{eq:balakrishnan_quadrature}. To do so, we pre-calculate $c_\text{const}$ from \eqref{eq:coefficients}, prepare a quadrature sum accumulator, and spawn an OpenMP parallel loop to calculate quadrature steps $l \in [-K^-, K^+]$. 
	
	\subsubsection*{The loop}
	For each quadrature step $l$ we compute \eqref{eq:coefficients}, and decompose the graph $\Gamma$ into individual edges. The resulting local operator, following \eqref{eq:spde_approx_fem} is
	\begin{equation} \label{eq:local_operator}
		A^{(l,e)} = c_{I}^{(l)} M^{(e)} + c_{L}^{(l)} S^{(e)}, \quad e \in E,
	\end{equation}
	or in block notation
	\begin{equation} \label{eq:spde_local_quadrature_block}
		\begin{bmatrix} A_{II}^{(l,e)} & A_{I\Gamma}^{(l,e)} \\ A_{\Gamma I}^{(l,e)} & A_{\Gamma\Gamma}^{(l,e)} \end{bmatrix} = 
		c_{I}^{(l)}
		\begin{bmatrix} M_{II}^{(e)} & M_{I\Gamma}^{(e)} \\ M_{\Gamma I}^{(e)} & M_{\Gamma\Gamma}^{(e)} \end{bmatrix}
		+ c_{L}^{(l)} 
		\begin{bmatrix} S_{II}^{(e)} & S_{I\Gamma}^{(e)} \\ S_{\Gamma I}^{(e)} & S_{\Gamma\Gamma}^{(e)} \end{bmatrix}
		\quad e \in E,
	\end{equation} 
	where the $\mathscr{C}_e (\cdot)$ operator returns the submatrix of its argument relevant to edge $e$. For each edge $e \in E$, the implementation avoids the overhead of sparse matrix formats. Instead, it allocates 1D arrays to store the tridiagonal components as $$d_{II}^{(l,e)} = \begin{bmatrix} A_{II,i,i}^{(l,e)} \end{bmatrix} \text{ for all } i
	\quad \quad
	d_{II,\up}^{(l,e)}, d_{II,\lo}^{(l,e)} = \begin{bmatrix} A_{II,i,j}^{(l,e)} \end{bmatrix}, \text{ for } j = i\pm1 .
	$$
	Since the internal nodes only connect to the boundary vertices at the ends of an edge, we store $A_{I \Gamma}^{(l,e)} = \left(A_{\Gamma I}^{(l,e)} \right)^T$, as a pair of scalar values representing the connections to $v_\text{left}$ and $v_\text{right}$.
	Next we utilize the Thomas Algorithm to obtain the factorization of the interior operator $A_{II}^{(l,e)} =LU$ with $L=\operatorname{diag}(d_{II,\lo}^{\prime (l,e)} \, , \, -1)$ and $U=\operatorname{diag}(d_{II}^{\prime (l,e)})+\operatorname{diag}(d_{II,\up}^{(l,e)}\, , \, +1)$, where
	\begin{align} \label{eq:thomas_coefficients}
		d_{II,1}^{\prime (l,e)}     & = d_{II,1}^{(l,e)},                                                                                     \\
		d_{II,\lo,i}^{\prime (l,e)} & = \frac{d_{II,\lo,i-1}^{(l,e)}}{d_{II,i-1}^{\prime (l,e)}}, \quad i \geq 2,                    \\
		d_{II,i}^{\prime (l,e)}     & = d_{II,i}^{(l,e)} - d_{II,\lo,i}^{\prime (l,e)} \cdot d_{II,\up,i-1}^{(l,e)}, \quad i \geq 2.
	\end{align}
	Analogously, we precompute $d^{(l,e)}, d_\lo^{(l,e)}$ and $d_\up^{(l,e)}$ from the global operator $A^{(l,e)}$.
	\subsubsection*{Boundary Condition and Initial Offset}
	By substituting \eqref{eq:spde_local_quadrature_block} into \eqref{eq:spde_quadrature_block} we obtain
	\begin{equation}  \label{eq:spde_edge_block}
		\begin{bmatrix} A_{II}^{(l,e)} & A_{I\Gamma}^{(l,e)} \\ A_{\Gamma I}^{(l,e)} & A_{\Gamma\Gamma}^{(l,e)} \end{bmatrix} \begin{bmatrix} u_I^{(l,e)} \\ u_\Gamma^{(l,e)} \end{bmatrix} = \begin{bmatrix} W_I^{(l,e)} \\ W_\Gamma^{(l,e)} \end{bmatrix}.
	\end{equation}
	In pursuit of reducing the global system to the topological nodes (i.e. to \emph{hide} the internal vertices) we express $u_I^{(l,e)}$ from the first block row of the matrix equations and substitute it into the second, resulting in the reduced local system:
	\begin{equation} \label{eq:local_reduced_schur}
		\underbrace{\left( A_{\Gamma\Gamma}^{(l,e)} - A_{\Gamma I}^{(l,e)} (A_{II}^{(l,e)})^{-1} A_{I\Gamma}^{(l,e)} \right)}_{S^{(l,e)}} u_\Gamma^{(l,e)} = \underbrace{
			W_\Gamma^{(l,e)} - A_{\Gamma I}^{(l,e)} 
			\overbrace{
				(A_{II}^{(l,e)})^{-1} W_I^{(l,e)}
			}^{w^{(l,e)}}
		}_{g_\Gamma^{(l,e)}}.
	\end{equation}
	We obtain $w^{(l,e)}$ as a solution to $A_{II}^{(l,e)} \cdot w^{(l,e)} =  W_I^{(l,e)}$ utilizing the pre-computed LU-factorization \eqref{eq:thomas_coefficients} $L \left( U w^{(l,e)} \right) = W_I^{(l,e)}$ in two steps:
	\begin{enumerate}
		\item Forwards Substitution ($Ly = W_I^{(l,e)}$)
		\begin{equation}
			y_1 = W_{I,1}^{(l,e)} \quad \text{then}\quad
			y_i = W_{I,i}^{(l,e)} - d_{II,\lo,i}^{\prime (l,e)} \cdot y_{i-1} \quad i = 2, \dots, n
		\end{equation}
		\item Backwards substitution ($U w^{(l,e)} = y$)
		\begin{equation}
			w^{(l,e)}_n = \frac{y_n}{d_{II,n}^{\prime (l,e)}} \quad \text{then} \quad
			w^{(l,e)}_i =
			\frac{y_i - d_{II,\up,i}^{(l,e)} \cdot w^{(l,e)}_{i+1}}
			{d_{II,i}^{\prime (l,e)}} \quad i= 1, \dots , (n-1)
		\end{equation}
	\end{enumerate}
	Since $A_{\Gamma I}^{(l,e)} \in \mathbb{R}^{2 \times (N_e - 2)}$ is extremely sparse with only two nonzero boundary elements $\left(A_{\Gamma I}^{(l,e)} \right)_{1,1}$ and $\left(A_{\Gamma I}^{(l,e)} \right)_{2,(N_e - 2)}$ we have that
	$$ g_\Gamma^{(l,e)} =
	\begin{bmatrix} 
		W_{\Gamma,1}^{(l,e)} - \left(A_{\Gamma I}^{(l,e)} \right)_{1,1} \cdot w^{(l,e)}_1                  \\ 
		W_{\Gamma,2}^{(l,e)} - \left(A_{\Gamma I}^{(l,e)} \right)_{2,(N_e - 2)} \cdot w^{(l,e)}_{(N_e -2)} 
	\end{bmatrix}. $$
	We further have that $g_\Gamma^{(l)} = \sum_{e \in E} \mathscr{S} \left( g_\Gamma^{(l,e)} \right)$ where the $\mathscr{S}(\cdot)$ operator expands its argument to a $|V| \times 1$ vector, routing the two local boundary values to their correct global vertex indices with zero fill-in. 
	\subsection{Phase 4: The Inside Loop (PCG solver)}
	Inside the PCG loop we are solving $S^{(l)} u_\Gamma^{(l)} = g_\Gamma^{(l)}$. We introduce the goal function $J(u_\Gamma^{(l)})$ we aim to minimze such that $\nabla J(u_\Gamma^{(l)}) = S^{(l)} u_\Gamma^{(l)} - g_\Gamma^{(l)}$. However, explicitly assembling $S^{(l)}$ would break the domain decomposition framework. We  demonstrate that we can construct an $\nabla J_2(u_\Gamma^{(l)})$, which yields the exact same numerical vector as $\nabla J_1$ but never requires the formation of $S^{(l)}$. The global matrix is the sum of all local Schur complements, expanded to the global dimensions $|V| \times |V|$ by padding them with zeros:
	\begin{equation} \label{eq:goal_func_1}
		\nabla J_1(u_\Gamma^{(l)}) = \left( \sum_{e \in E} \mathscr{S} \left( S^{(l,e)} \right) \right) u_\Gamma^{(l)} - g_\Gamma^{(l)} = \sum_{e \in E} \mathscr{S} \left( S^{(l,e)} u_\Gamma^{(l,e)} \right) - g_\Gamma^{(l)}.
	\end{equation}
	Next, we substitute the definition of the local Schur complement $S^{(l,e)}$ from \eqref{eq:local_reduced_schur} to evaluate the local matrix-vector product $y^{(l,e)} = S^{(l,e)} u_\Gamma^{(l,e)}$:
	\begin{equation} \label{eq:matrix_free_local}
		y^{(l,e)} = A_{\Gamma\Gamma}^{(l,e)} u_\Gamma^{(l,e)} - A_{\Gamma I}^{(l,e)} \underbrace{\left[ (A_{II}^{(l,e)})^{-1} (A_{I\Gamma}^{(l,e)} u_\Gamma^{(l,e)}) \right]}_{z^{(l,e)}}.
	\end{equation}
	\subsubsection*{The Dirichlet Step}
	This substitution allows us to define our algorithmic gradient $\nabla J_2$ entirely in terms of local, computationally inexpensive operations. For the inner term $z^{(l,e)}$, rather than computing a dense inverse matrix $(A_{II}^{(l,e)})^{-1}$, we solve the equivalent tridiagonal system $A_{II}^{(l,e)} z^{(l,e)} = A_{I\Gamma}^{(l,e)} u_\Gamma^{(l,e)}$ using the cached Thomas factorization from the setup phase. 
	Finally, we substitute back into \eqref{eq:goal_func_1} to obtain
	\begin{equation} \label{eq:goalfunction}
		\nabla J_2(u_\Gamma^{(l)}) \equiv \nabla J_1(u_\Gamma^{(l)}) \equiv \nabla J(u_\Gamma^{(l)}) = \sum_{e \in E} \mathscr{S} \left( y^{(l,e)} \right) - g_\Gamma^{(l)}.
	\end{equation}
	\subsubsection*{The Neumann Step}
	Having computed the global topological vertex flux residual vector $r_\Gamma^{(l)} = \nabla J_2(u_\Gamma^{(l)})$ for a given $u_\Gamma^{(l)}$, we distribute it across edges inversely proportional to te degree $d_v$ of each vertex $v$.
	\begin{equation}\label{eq:distribute_dv}
		\tilde{r}_\Gamma^{(l,e)} = \mathscr{C}_e \left( r_\Gamma^{(l)} \operatorname{diag} \left( [\dots \frac{1}{d_v} \dots] \right) \right) \quad v \in V.
	\end{equation}
	We consider
	\begin{equation} \label{eq:neumann_block}
		\underbrace{\begin{bmatrix} A_{II}^{(l,e)} & A_{I\Gamma}^{(l,e)} \\ A_{\Gamma I}^{(l,e)} & A_{\Gamma\Gamma}^{(l,e)} \end{bmatrix}}_\text{Not tridiagonal} \begin{bmatrix} \delta_I^{(l,e)} \\ \delta_\Gamma^{(l,e)} \end{bmatrix} = \begin{bmatrix} 0 \\ \tilde{r}^{(l,e)} \end{bmatrix}.
	\end{equation}
	To rearrange the above system of equations into a tridiagonal one, we must apply the permutation operator
	$$
	P = \begin{bmatrix} 
		0_{1 \times (N-2)} & 1 & 0 \\ 
		I_{(N-2) \times (N-2)} & 0_{(N-2) \times 1} & 0_{(N-2) \times 1} \\ 
		0_{1 \times (N-2)} & 0 & 1 
	\end{bmatrix}
	$$
	to \eqref{eq:neumann_block} so as to obtain
	\begin{equation*}
		\label{eq:neumann_block_tridiag}
		\underbrace{P A^{(l,e)} P^T}_\text{Tridiagonal} 
		\begin{bmatrix} 
			\delta_{\Gamma,\text{left}}^{(l,e)}   \\
			\delta_I^{(l,e)}                      \\ 
			\delta_{\Gamma,\text{right}}^{(l,e)} 
		\end{bmatrix}=
		\begin{bmatrix} 
			\tilde{r}_\text{left}^{(l,e)}  \\ 
			0_{1 \times (N-2)}             \\ 
			\tilde{r}_\text{right}^{(l,e)} 
		\end{bmatrix},
	\end{equation*}
	which can efficiently be solved by the Thomas Algorithm. For computational efficiency, we never explicitly build the permutation matrix. Instead, we directly populate 1D vectors representing the diagonals of the transformed matrix. After computing the solution vector via the Thomas algorithm, we extract the boundary displacements $\delta_{\Gamma,1}^{(l,e)}$ and $\delta_{\Gamma,2}^{(l,e)}$, and omit the internal components. Then, we consider the averaged flux residual induced displacements
	\begin{equation} \label{eq:induced_displacements}
		w_\Gamma^{(l)} = \sum_{e \in E} \mathscr{C}_e \left( \begin{bmatrix} 1/d_\text{left} & 0 \\ 0 & 1/d_\text{right} \end{bmatrix} \begin{bmatrix} \delta_{\Gamma,\text{left}}^{(l,e)} \\ \delta_{\Gamma,\text{right}}^{(l,e)} \end{bmatrix} \right)
	\end{equation}
	in each topological node as the preconditioned gradient contributing to the PCG solver's subsequent solution estimate. The PCG loop continues this loop until a solution $ J(u_\Gamma^{(l)}) \approx 0_{|V| \times 1}$ is reached.
	
	Note that the Neumann Step is equivalent to applying the preconditioner $D_\Gamma^{-1} \left( \sum_{e \in E} (S^{(l,e)})^{-1} \right) D_\Gamma^{-1}$ as in \cite[Section 3.2.1]{Kovacs_2025}, where the diagonal elements of $D_\Gamma$ are $d_v$ for $v \in \Gamma$. Distributing the residual inversely proportional to $d_v$ in \eqref{eq:distribute_dv} applies the first diagonal scaling $D_\Gamma^{-1}$. Solving the local block system \eqref{eq:neumann_block} computes the effect of the local inverse Schur complement $(S^{(l,e)})^{-1}$ as $\delta_\Gamma^{(l,e)} = (S^{(l,e)})^{-1} \tilde{r}^{(l,e)}$. Finally, averaging the boundary displacements in \eqref{eq:induced_displacements} applies the outer $D_\Gamma^{-1}$ scaling.
	\subsection{Phase 5: Post-Processing}
	\subsubsection*{Recovery of Internal Nodal Displacements}
	Because the PCG solver operates on the topological vertices of the metric graph, we must recover the corresponding internal nodal values $u_I^{(l,e)}$ for every edge $e \in E$.
	We rearrange the first block row of the original un-reduced local system from \eqref{eq:spde_edge_block} as
	\begin{equation} \label{eq:internal_recovery}
		A_{II}^{(l,e)} u_I^{(l,e)} = W_I^{(l,e)} -
		A_{I\Gamma}^{(l,e)} u_\Gamma^{(l,e)}.
	\end{equation}
	Since $A_{I\Gamma}^{(l,e)} u_\Gamma^{(l,e)}$ only has two nonzero elements in the first and last positions, we directly modify the same positions of $W_I^{(l,e)}$ to obtain the right-hand-side.
	Next, we utilize the prefactorization of $A_{II}^{(l,e)}$ and calculate $u_I^{(l,e)}$ using the Thomas algorithm.
	\subsubsection*{Quadrature Accumulation and Global Assembly}
	With the topological vertex solutions $u_\Gamma^{(l)}$ and the recovered internal solution vectors $ u_I^{(l,e)}$, we form 
	\begin{equation*}	
		u^{(l)} = \begin{bmatrix} u_{I}^{(l)} ~ u_\Gamma^{(l)} \end{bmatrix}^T = 
		\begin{bmatrix}
			u_{I}^{(l, e_1)} ~    
			\dots ~                
			u_{I}^{(l, e_{|E|})} ~ 
			u_\Gamma^{(l)}         
		\end{bmatrix}^T .
	\end{equation*}
	Finally, we accumulate our results over all quadrature steps as per \eqref{eq:spde_approx_fem}: $u = \sum_{l = -K^-}^{K^+} u^{(l)}$.
	\section{Numerical Experiments} \label{sec:numerical}
	The numerical framework was implemented in C++, utilizing the \texttt{Eigen3} library for all linear algebra operations. To parallelize the independent spatial equations generated by the Balakrishnan quadrature, we employed OpenMP. 
	All benchmarks were orchestrated using the Slurm Workload Manager on a high-performance computing cluster. To guarantee a reproducible execution environment, the application was deployed within an Ubuntu-based Apptainer container. Furthermore, to optimize the heavy dynamic memory allocation overhead inherent to assembling large sparse matrices, \texttt{glibc} allocator was overridden with Google's thread-caching alternative \texttt{tcmalloc}. Finally, strict thread affinity (\texttt{OMP\_PROC\_BIND=spread} and \texttt{OMP\_PLACES=threads}) was enforced at runtime to prevent operating system thread-migration overhead.
	To ensure less machine-dependent metrics, performance measurement were run strictly single-threaded, and computational cost was measured via hardware instruction counts (\texttt{PAPI\_TOT\_INS}) using the Performance Application Programming Interface (PAPI) library. Memory consumption was quantified by tracking the peak Resident Set Size (RSS), measured directly from the Linux \texttt{/proc/self/statm} pseudo-filesystem during the most memory-intensive lifecycle phases of the solver.
	\subsection{Strong error}
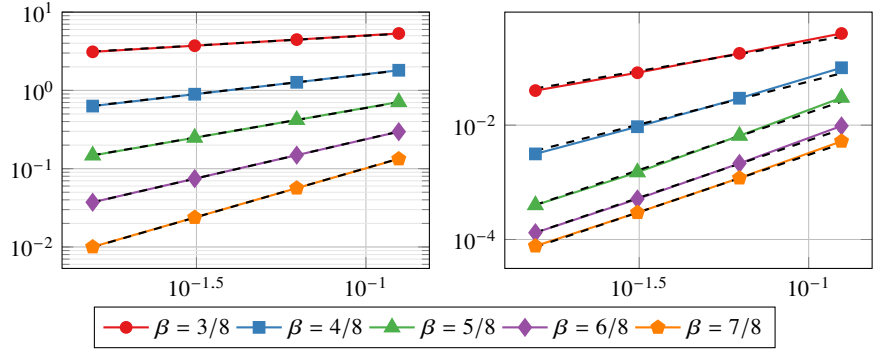
\begin{figure}[htbp]
		\centering	
		\begin{subfigure}[b]{0.5\textwidth}
			\centering
			\begin{tikzpicture}
				\begin{axis}[
					width=1.1*\linewidth, height=0.85\linewidth,
					xmode=log, ymode=log,
					grid=both, grid style={line width=.1pt, draw=gray!20},
					major grid style={line width=.2pt,draw=gray!50},
					legend to name=sharedlegend,
					legend style={
						legend columns=5, 
						cells={anchor=west}, font=\footnotesize
					},
					x dir=normal
					]
					
					\addplot[color=color1, mark=*, thick] table[x=h, y=err375] {\datatableStrong};
					\addlegendentry{$\beta=3/8$}
					\addplot[color=black, dashed, thick, domain=0.015625:0.125000, forget plot] {8.881 * x^(0.25)};
					
					\addplot[color=color2, mark=square*, thick] table[x=h, y=err500] {\datatableStrong};
					\addlegendentry{$\beta=4/8$}
					\addplot[color=black, dashed, thick, domain=0.015625:0.125000, forget plot] {5.071 * x^(0.5)};
					
					\addplot[color=color3, mark=triangle*, thick, mark size=3pt] table[x=h, y=err625] {\datatableStrong};
					\addlegendentry{$\beta=5/8$}
					\addplot[color=black, dashed, thick, domain=0.015625:0.125000, forget plot] {3.353 * x^(0.75)};
					
					\addplot[color=color4, mark=diamond*, thick, mark size=3pt] table[x=h, y=err750] {\datatableStrong};
					\addlegendentry{$\beta=6/8$}
					\addplot[color=black, dashed, thick, domain=0.015625:0.125000, forget plot] {2.381 * x^(1.0)};
					
					\addplot[color=color5, mark=pentagon*, thick, mark size=2.5pt] table[x=h, y=err875] {\datatableStrong};
					\addlegendentry{$\beta=7/8$}
					\addplot[color=black, dashed, thick, domain=0.015625:0.125000, forget plot] {1.805 * x^(1.25)};
					
				\end{axis}
			\end{tikzpicture}
		\end{subfigure}\hfill
		\begin{subfigure}[b]{0.5\textwidth}
			\centering
			\begin{tikzpicture}
				\begin{axis}[
					width=1.1*\linewidth, height=0.85\linewidth,
					xmode=log, ymode=log,
					grid=both, grid style={line width=.1pt, draw=gray!20},
					major grid style={line width=.2pt,draw=gray!50},
					x dir=normal
					]
					
					\addplot[color=color1, mark=*, thick] table[x=h, y=cov375] {\datatableCov};
					\addplot[color=black, dashed, thick, domain=0.015625:0.125000, forget plot] {2.806 * x^(1.0)};
					
					\addplot[color=color2, mark=square*, thick] table[x=h, y=cov500] {\datatableCov};
					\addplot[color=black, dashed, thick, domain=0.015625:0.125000, forget plot] {1.828 * x^(1.5)};
					
					\addplot[color=color3, mark=triangle*, thick, mark size=3pt] table[x=h, y=cov625] {\datatableCov};
					\addplot[color=black, dashed, thick, domain=0.015625:0.125000, forget plot] {1.648 * x^(2.0)};
					
					\addplot[color=color4, mark=diamond*, thick, mark size=3pt] table[x=h, y=cov750] {\datatableCov};
					\addplot[color=black, dashed, thick, domain=0.015625:0.125000, forget plot] {0.540 * x^(2.0)};
					
					\addplot[color=color5, mark=pentagon*, thick, mark size=2.5pt] table[x=h, y=cov875] {\datatableCov};
					\addplot[color=black, dashed, thick, domain=0.015625:0.125000, forget plot] {0.302 * x^(2.0)};
					
				\end{axis}
			\end{tikzpicture}
		\end{subfigure}	
		\vspace{0.1cm} 
		\ref{sharedlegend}
		\vspace{0.1cm}
		%
		\caption{Observed strong error (left) and covariance error (right). The horizontal axes represent mesh size $h$. Dashed lines indicate the exact theoretical rates.}
		\label{fig:convergence_results}
	\end{figure}

	To validate whether the lumped-mass-based method breaks convergence metrics, we set out to reproduce the experiments from \cite[Chapter 7]{Bolin_2023}, but with mass lumping. According to them, the theoretical rate of strong convergence is $2 \beta - \frac{1}{2}$ for $\beta \in (\frac{1}{4},1)$. We considered the fractional SPDE \eqref{eq:fractional_spde} with $\kappa =1$, quadrature step size
	$k = \frac{-1}{\beta \ln h}$, quadrature limits $K^- = \left\lceil \frac{\pi^2}{4k^2\beta} \right\rceil$, $K^+ = \left\lceil \frac{\pi^2}{4k^2(1-\beta)} \right\rceil$ and $\beta = \frac{n}{8}$ for $n=3,4,5,6$. Every edge in the graph was discretized into equally sized segments such that $h_\text{max} = 2^{-\ell}$. For the overkill (or reference) solution $u_\text{ok}$ we set $h_{\text{ok},\text{max}} = 2^{-16}$. Closely following \cite{Bolin_2023} we project the noise generated on the overkill mesh down, we solve the coarse problem, upsample the result, and measure the $L^2(\Gamma)$ norm difference as compared to the reference solution. For each $\beta$ value this process is repeated 100 times with different noise realizations on the fine mesh. Then, we take the average squared $L^2$ errors $ \text{err} = \frac{1}{100}\sum_{i=1}^{100} (u_\text{ok}^{(i)} - u_\text{coarse}^{(i)})^T {M}_\text{ok} (u_\text{ok}^{(i)} - u_\text{coarse}^{(i)})$. For each $\beta$ then, the estimated rate of convergence $r$ is expressed from the linear regression $\ln \text{err} = c + r \ln h$ as depicted in the left panel of \figref{fig:convergence_results} and in \tabref{tab:convergence_rates}.
	\subsection{Covariance error}
	Similarly, \cite{Bolin_2023} shows that the theoretical rate of convergence of the $L^2(\Gamma \times \Gamma)$ covariance error is $\min (4 \beta - \frac{1}{2},2)$ for $\beta \in (\frac{1}{4},1)$ for $k = \frac{-1}{\beta \ln h}$. For the same overkill and coarse mesh resolutions and $\beta$ values, we considered piecewise constant covariance functions on the overkill mesh $h=2^{-8}$. The corresponding results are depicted in the right panel of \figref{fig:convergence_results} and in \tabref{tab:convergence_rates}.
	\subsection{Performance comparison}
\begin{figure}[htbp]
		\centering
		\begin{subfigure}[b]{0.5\textwidth}
			\centering
			\begin{tikzpicture}
				\begin{axis}[
					title={Total Instruction Count},
					xmode=log,
					ymode=log,
					grid=both,
					minor grid style={gray!25},
					major grid style={gray!50},
					width=1.1*\linewidth,
					height=0.9\linewidth, 
					mark size=2.5pt,
					legend to name=sharedperflegend,
					legend style={
						legend columns=4, 
						cells={anchor=west}, font=\footnotesize
					}
					]
					
					\addplot[color=blue, mark=square, thick] 
					table[x=N, y=Instr_L] {\datatableFull};
					\addlegendentry{Lumped (Full)}
					
					\addplot[color=red, mark=square, thick] 
					table[x=N, y=Instr_C] {\datatableFull};
					\addlegendentry{Consistent (Full)}
					
					\addplot[color=blue!70, mark=x, mark size=3pt, mark options={solid}, thick, dashed] 
					table[x=N, y=Instr_L] {\datatableNoise};
					\addlegendentry{Lumped (Noise)}
					
					\addplot[color=red!70, mark=x, mark size=3pt, mark options={solid}, thick, dashed] 
					table[x=N, y=Instr_C] {\datatableNoise};
					\addlegendentry{Consistent (Noise)}
					
				\end{axis}
			\end{tikzpicture}
		\end{subfigure}\hfill
		\begin{subfigure}[b]{0.5\textwidth}
			\centering
			\begin{tikzpicture}
				\begin{axis}[
					title={Peak Memory Consumption [KiB]},
					xmode=log,
					ymode=log,
					grid=both,
					minor grid style={gray!25},
					major grid style={gray!50},
					width=1.1*\linewidth,
					height=0.9\linewidth,
					mark size=2.5pt
					]
					
					\addplot[color=blue, mark=square, thick] 
					table[x=N, y expr=\thisrow{Mem_L}] {\datatableFull};
					
					\addplot[color=red, mark=square, thick] 
					table[x=N, y expr=\thisrow{Mem_C}] {\datatableFull};
					
					\addplot[color=blue!70, mark=x, mark size=3pt, mark options={solid}, thick,dashed] 
					table[x=N, y expr=\thisrow{Mem_L}] {\datatableNoise};
					
					\addplot[color=red!70, mark=x, mark size=3pt, mark options={solid}, thick,dashed] 
					table[x=N, y expr=\thisrow{Mem_C}] {\datatableNoise};
					
				\end{axis}
			\end{tikzpicture}
		\end{subfigure}
		\vspace{0.1cm} 
		\ref{sharedperflegend}
		\vspace{0.1cm}
		%
		\caption{Performance comparison the GRF generation pipeline between Lumped and Consistent mass formulations. The horizontal axes represent the total number of topological vertices in the graph.}
		\label{fig:performance_scaling_merged}
	\end{figure}
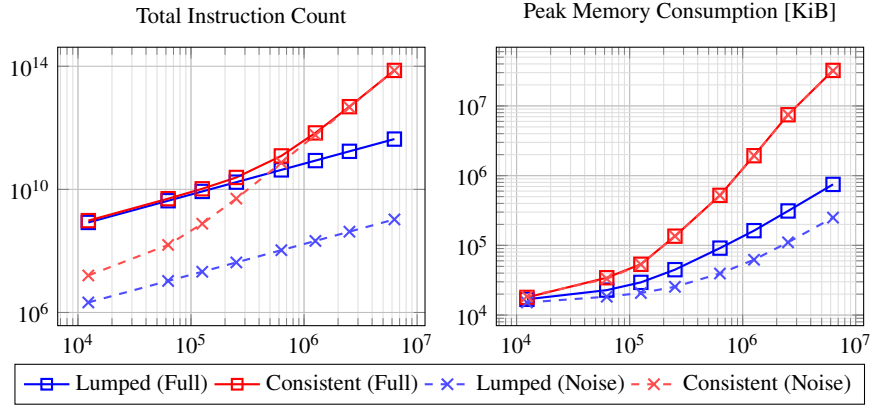

	Finally, we conducted a performance scaling analysis to quantify the computational advantages of the lumped mass formulation over the consistent mass approach. We generated Barabási-Albert scale-free metric graphs with attachment parameter $m=2$ and with increasing topological complexity, ranging from 100 to 50,000 vertices, while keeping the mesh resolution fixed at $h = 2^{-6}$. The computational cost was measured via hardware instruction counts using the PAPI library, whereas for memory consumption we considered the peak Resident Set Memory metric of the running process. Results, depicted in \figref{fig:performance_scaling_merged}, are evaluated across two scenarios: isolated spatial white noise generation, and the full GRF generation including the Balakrishnan quadrature and PCG solver.
	\pgfplotstableread{
		h           err375      err500      err625      err750      err875
		0.125000    5.335224    1.807369    0.710203    0.296709    0.133258
		0.062500    4.455675    1.268635    0.419348    0.148768    0.056451
		0.031250    3.721420    0.893495    0.248365    0.074549    0.023789
		0.015625    3.107954    0.630631    0.147539    0.037264    0.010011
	}\datatableStrong
	\pgfplotstableread{
		h           cov375      cov500      cov625      cov750      cov875
		0.125000    0.40125300  0.10082900  0.03016750  0.00969306  0.00522219
		0.062500    0.18013200  0.02957370  0.00651719  0.00214048  0.00118818
		0.031250    0.08249470  0.00938936  0.00149840  0.00051134  0.00029512
		0.015625    0.04023120  0.00314021  0.00039892  0.00013150  0.00007682
	}\datatableCov
	\begin{table}[H]
		\centering
		\setlength{\tabcolsep}{7pt}
		\begin{tabular}{l c c c c c}
			\toprule
			$\beta$           & $3/8$       & $4/8$       & $5/8$       & $6/8$       & $7/8$       \\
			\midrule
			Strong rate & 0.26 (0.25) & 0.51 (0.50) & 0.76 (0.75) & 1.00 (1.00) & 1.25 (1.25) \\
			Covariance rate   & 1.07 (1.0)  & 1.59 (1.5)  & 2.00 (2.0)  & 2.01 (2.0)  & 1.98 (2.0)  \\
			\bottomrule
		\end{tabular}
		\caption{Empirical convergence rates. Theoretical rates are given in parentheses.}
		\label{tab:convergence_rates}
	\end{table}	
	\pgfplotstableread{
		Graph                     Task         N_Graph  h              N        Instr_L         Iters_L  Mem_L      Instr_C         Iters_C  Mem_C      Speedup
		barabasi_albert_100       noise_only   100      1.562500e-02   12348    2116264         0        15183      15929913        0        18512      7.53
		barabasi_albert_500       noise_only   500      1.562500e-02   62748    10586941        0        18359      156738266       0        33194      14.80
		barabasi_albert_1000      noise_only   1000     1.562500e-02   125748   21053787        0        20739      765501694       0        53182      36.36
		barabasi_albert_2000      noise_only   2000     1.562500e-02   251748   42247641        0        25457      5003188030      0        135725     118.43
		barabasi_albert_5000      noise_only   5000     1.562500e-02   629748   105559643       0        39476      73784101998     0        520391     698.98
		barabasi_albert_10000     noise_only   10000    1.562500e-02   1259748  211144959       0        62182      586339217535    0        1931610    2776.95
		barabasi_albert_20000     noise_only   20000    1.562500e-02   2519748  422335767       0        110518     4671046939967   0        7473586    11060.03
		barabasi_albert_50000     noise_only   50000    1.562500e-02   6299748  1055734441      0        250008     72752990051284  0        31923868   68912.21
	}\datatableNoise	
	\pgfplotstableread{
		Graph                     Task         N_Graph  h              N        Instr_L         Iters_L  Mem_L      Instr_C         Iters_C  Mem_C      Speedup
		barabasi_albert_100       full_nrf     100      1.562500e-02   12348    840973767       289      16845      965435287       289      17943      1.15
		barabasi_albert_500       full_nrf     500      1.562500e-02   62748    4272453675      290      22798      4981034452      290      34352      1.17
		barabasi_albert_1000      full_nrf     1000     1.562500e-02   125748   8559598704      290      29466      10433153541     290      53574      1.22
		barabasi_albert_2000      full_nrf     2000     1.562500e-02   251748   17148361896     290      44810      24369933376     290      135477     1.42
		barabasi_albert_5000      full_nrf     5000     1.562500e-02   629748   42949809481     290      91316      122283211848    290      521666     2.85
		barabasi_albert_10000     full_nrf     10000    1.562500e-02   1259748  86024571167     290      162509     683464726693    290      1932120    7.94
		barabasi_albert_20000     full_nrf     20000    1.562500e-02   2519748  172235423490    290      311318     4865487119677   290      7472514    28.25
		barabasi_albert_50000     full_nrf     50000    1.562500e-02   6299748  430610178036    290      749062     73239121672119  290      32296509   170.08
	}\datatableFull
	\section{Conclusion} \label{sec:conclusion}
	In this paper, we presented a numerical framework for generating Gaussian Random Fields on metric graphs. We coupled a Neumann-Neumann domain decomposition method with finite element mass matrix lumping, mitigating the severe computational and memory bottlenecks inherent to the standard finite element SPDE approach. Our numerical experiments confirm that mass lumping does not degrade the convergence rates for strong error and covariance error, but perfectly mirror the exact theoretical rates established for the consistent mass formulation.
	Our numerical results show that as the graph size increases, the dense fill-in of Cholesky factorization exhibits a non-linear computational explosion -- both in peak memory consumption and operation count -- completely dominating the GRF generation pipeline. Conversely, mass lumping diagonalizes the noise projection, reducing its computational complexity to strictly $\mathcal{O}(N)$ operations, resulting in a linear relationship between the computational cost of full GRF generation pipeline and topological node count. We have furthermore demonstrated that algorithmic efficiency extends directly to the memory footprint, resulting in memory savings beyond an order of magnitude.
	\section*{Acknowledgment} M. Kovács acknowledges the support of the Hungarian
	National Research, Development and Innovation Office (NKFIH) through Grant no. K-145934.
	\bibliographystyle{splncs04}	
	\bibliography{bibliography}

@Article{Arioli_2017,
  author    = {Arioli, Mario and Benzi, Michele},
  journal   = {IMA Journal of Numerical Analysis},
  title     = {A finite element method for quantum graphs},
  year      = {2017},
  issn      = {1464-3642},
  month     = jun,
  number    = {3},
  pages     = {1119--1163},
  volume    = {38},
  doi       = {10.1093/imanum/drx029},
  publisher = {Oxford University Press (OUP)},
}

@Article{Lindgren_2011,
  author    = {Lindgren, Finn and Rue, Håvard and Lindström, Johan},
  journal   = {Journal of the Royal Statistical Society Series B: Statistical Methodology},
  title     = {An Explicit Link between Gaussian Fields and Gaussian Markov Random Fields: The Stochastic Partial Differential Equation Approach},
  year      = {2011},
  issn      = {1467-9868},
  month     = aug,
  number    = {4},
  pages     = {423--498},
  volume    = {73},
  doi       = {10.1111/j.1467-9868.2011.00777.x},
  publisher = {Oxford University Press (OUP)},
}

@Article{Bolin_2024,
  author    = {Bolin, David and Simas, Alexandre B. and Wallin, Jonas},
  journal   = {Bernoulli},
  title     = {Gaussian Whittle–Matérn fields on metric graphs},
  year      = {2024},
  issn      = {1350-7265},
  month     = may,
  number    = {2},
  volume    = {30},
  doi       = {10.3150/23-bej1647},
  publisher = {Bernoulli Society for Mathematical Statistics and Probability},
}

@Article{Bolin_2023,
  author    = {Bolin, David and Kovács, Mihály and Kumar, Vivek and Simas, Alexandre},
  journal   = {Mathematics of Computation},
  title     = {Regularity and numerical approximation of fractional elliptic differential equations on compact metric graphs},
  year      = {2023},
  issn      = {1088-6842},
  month     = dec,
  number    = {349},
  pages     = {2439--2472},
  volume    = {93},
  doi       = {10.1090/mcom/3929},
  publisher = {American Mathematical Society (AMS)},
}

@Article{Kovacs_2025,
  author    = {Kovács, Mihály and Vághy, Mihály},
  journal   = {BIT Numerical Mathematics},
  title     = {Neumann-Neumann type domain decomposition of elliptic problems on metric graphs},
  year      = {2025},
  issn      = {1572-9125},
  month     = may,
  number    = {2},
  volume    = {65},
  doi       = {10.1007/s10543-025-01067-8},
  publisher = {Springer Science and Business Media LLC},
}

@Book{Toselli_2005,
  author    = {Toselli, Andrea and Widlund, Olof B.},
  publisher = {Springer Berlin Heidelberg},
  title     = {Domain Decomposition Methods — Algorithms and Theory},
  year      = {2005},
  isbn      = {9783540266624},
  doi       = {10.1007/b137868},
  issn      = {2198-3712},
  journal   = {Springer Series in Computational Mathematics},
}

@Book{Berkolaiko_2012,
  author    = {Berkolaiko, Gregory and Kuchment, Peter},
  publisher = {American Mathematical Society},
  title     = {Introduction to Quantum Graphs},
  year      = {2012},
  isbn      = {9780821894552},
  month     = dec,
  doi       = {10.1090/surv/186},
  issn      = {2331-7159},
  journal   = {Mathematical Surveys and Monographs},
}

@Book{Davis_2006,
  author    = {Davis, Timothy A.},
  publisher = {Society for Industrial and Applied Mathematics},
  title     = {Direct Methods for Sparse Linear Systems},
  year      = {2006},
  isbn      = {9780898718881},
  month     = jan,
  doi       = {10.1137/1.9780898718881},
}

@Book{Saad_2003,
  author    = {Saad, Yousef},
  publisher = {Society for Industrial and Applied Mathematics},
  title     = {Iterative Methods for Sparse Linear Systems},
  year      = {2003},
  isbn      = {9780898718003},
  month     = jan,
  doi       = {10.1137/1.9780898718003},
}

@Book{Golub_2013,
  author    = {Golub, Gene H. and Loan, Charles F. Van},
  publisher = {Johns Hopkins Univ. Press},
  title     = {Matrix computations},
  year      = {2013},
  address   = {Baltimore, MD},
  edition   = {4. ed.},
  isbn      = {9781421407944},
  note      = {In association with the Department of Mathematical Sciences, The Johns Hopkins University.},
  series    = {Johns Hopkins studies in the mathematical sciences},
  pagetotal = {756},
  ppn_gvk   = {1957377690},
}

@Article{Balakrishnan_1960,
  author    = {Balakrishnan, A.},
  journal   = {Pacific Journal of Mathematics},
  title     = {Fractional powers of closed operators and the semigroups generated by them},
  year      = {1960},
  issn      = {0030-8730},
  month     = jun,
  number    = {2},
  pages     = {419--437},
  volume    = {10},
  doi       = {10.2140/pjm.1960.10.419},
  publisher = {Mathematical Sciences Publishers},
}

@Article{Bonito_2015,
  author    = {Bonito, Andrea and Pasciak, Joseph E.},
  journal   = {Mathematics of Computation},
  title     = {Numerical approximation of fractional powers of elliptic operators},
  year      = {2015},
  issn      = {1088-6842},
  month     = mar,
  number    = {295},
  pages     = {2083--2110},
  volume    = {84},
  doi       = {10.1090/s0025-5718-2015-02937-8},
  publisher = {American Mathematical Society (AMS)},
}

@Book{Quarteroni_1994,
  author    = {Quarteroni, Alfio and Valli, Alberto},
  publisher = {Springer Berlin Heidelberg},
  title     = {Numerical Approximation of Partial Differential Equations},
  year      = {1994},
  isbn      = {9783540852681},
  doi       = {10.1007/978-3-540-85268-1},
  issn      = {0179-3632},
  journal   = {Springer Series in Computational Mathematics},
}

@Article{Alexander_1983,
  author    = {Alexander, S.},
  journal   = {Physical Review B},
  title     = {Superconductivity of networks. A percolation approach to the effects of disorder},
  year      = {1983},
  issn      = {0163-1829},
  month     = feb,
  number    = {3},
  pages     = {1541--1557},
  volume    = {27},
  doi       = {10.1103/physrevb.27.1541},
  publisher = {American Physical Society (APS)},
}

@Article{Bouchaud_1987,
  author    = {Bouchaud, J. P and Lhuillier, C},
  journal   = {Europhysics Letters (EPL)},
  title     = {High-Field Behaviour of Liquid and Solid 3 He. A New Solid Phase?},
  year      = {1987},
  issn      = {1286-4854},
  month     = feb,
  number    = {4},
  pages     = {481--488},
  volume    = {3},
  doi       = {10.1209/0295-5075/3/4/015},
  publisher = {IOP Publishing},
}

@Article{Flesia_1989,
  author    = {Flesia, Cristina and Johnston, Robert and Kunz, Hervé},
  journal   = {Physical Review A},
  title     = {Localization of classical waves in a simple model},
  year      = {1989},
  issn      = {0556-2791},
  month     = oct,
  number    = {7},
  pages     = {4011--4018},
  volume    = {40},
  doi       = {10.1103/physreva.40.4011},
  publisher = {American Physical Society (APS)},
}

@Article{Kallianpur_1984,
  author    = {Kallianpur, G. and Wolpert, R.},
  journal   = {Applied Mathematics \& Optimization},
  title     = {Infinite dimensional stochastic differential equation models for spatially distributed neurons},
  year      = {1984},
  issn      = {1432-0606},
  month     = oct,
  number    = {1},
  pages     = {125--172},
  volume    = {12},
  doi       = {10.1007/bf01449039},
  publisher = {Springer Science and Business Media LLC},
}

@Article{Cho_2018,
  author    = {Cho, H. and Ayers, K. and de Pills, L. and Kuo, Y. and Park, J. and Radunskaya, A. and Rockne, R.},
  journal   = {Letters in Biomathematics},
  title     = {Modelling Acute Myeloid Leukaemia in a Continuum of Differentiation States},
  year      = {2018},
  issn      = {2373-7867},
  number    = {2},
  volume    = {5},
  doi       = {10.30707/lib5.2cho},
  publisher = {Illinois State University},
}

@Article{Avdonin_2023,
  author    = {Avdonin, Sergei and Edward, Julian and Leugering, Günter},
  journal   = {Evolution Equations and Control Theory},
  title     = {Controllability for the wave equation on graph with cycle and delta-prime vertex conditions},
  year      = {2023},
  issn      = {2163-2480},
  number    = {6},
  pages     = {1542--1558},
  volume    = {12},
  doi       = {10.3934/eect.2023025},
  publisher = {American Institute of Mathematical Sciences (AIMS)},
}

@Article{Avdonin_2019,
  author    = {Avdonin, Sergei and Zhao, Yuanyuan},
  journal   = {Applied Mathematics \& Optimization},
  title     = {Exact Controllability of the 1-D Wave Equation on Finite Metric Tree Graphs},
  year      = {2019},
  issn      = {1432-0606},
  month     = nov,
  number    = {3},
  pages     = {2303--2326},
  volume    = {83},
  doi       = {10.1007/s00245-019-09629-3},
  publisher = {Springer Science and Business Media LLC},
}

@Article{Mehandiratta_2021,
  author    = {Mehandiratta, Vaibhav and Mehra, Mani and Leugering, Gunter},
  journal   = {SIAM Journal on Control and Optimization},
  title     = {Optimal Control Problems Driven by Time-Fractional Diffusion Equations on Metric Graphs: Optimality System and Finite Difference Approximation},
  year      = {2021},
  issn      = {1095-7138},
  month     = jan,
  number    = {6},
  pages     = {4216--4242},
  volume    = {59},
  doi       = {10.1137/20m1340332},
  publisher = {Society for Industrial & Applied Mathematics (SIAM)},
}

@Article{Stoll_2021,
  author    = {Stoll, Martin and Winkler, Max},
  journal   = {ETNA - Electronic Transactions on Numerical Analysis},
  title     = {Optimal Dirichlet control of partial differential equations on networks},
  year      = {2021},
  issn      = {1068-9613},
  pages     = {392--419},
  volume    = {54},
  doi       = {10.1553/etna\_vol54s392},
  publisher = {Osterreichische Akademie der Wissenschaften},
}

@Article{AMR,
  author    = {Anderes, Ethan and Møller, Jesper and Rasmussen, Jakob G.},
  journal   = {The Annals of Statistics},
  title     = {Isotropic covariance functions on graphs and their edges},
  year      = {2020},
  issn      = {0090-5364},
  month     = Aug,
  number    = {4},
  volume    = {48},
  doi       = {10.1214/19-aos1896},
  publisher = {Institute of Mathematical Statistics},
}
\end{document}